\documentclass[10pt]{amsart}
\usepackage[all]{xy}
\usepackage{amsmath}
\usepackage{amsfonts}
\usepackage{amssymb}
\usepackage{amscd}
\usepackage{verbatim}
\newtheorem{theorem}{Theorem} [section]
\newtheorem{lemma} [theorem] {Lemma}
\newtheorem{proposition} [theorem] {Proposition}

\newtheorem{corollary} [theorem] {Corollary}
\newtheorem{definition} [theorem] {Definition}

\newtheorem{claim}[theorem] {Claim}
\begin{document}
\title{Semiglobal results for $\overline\partial$ on complex spaces with arbitrary singularities, Part II}
\author{Nils \O vrelid and Sophia Vassiliadou}
\thanks{The research of the second author is partially supported by NSF grant DMS-0712795}
\address{Dept. of Mathematics\\University of Oslo\\
P.B 1053 Blindern, Oslo, N-0316 NORWAY}
\address{Dept. of Mathematics\\Georgetown University\\
Washington, DC 20057 USA} \email{nilsov@math.uio.no,\;
sv46@georgetown.edu}
\date{\today}
\begin{abstract} \noindent We obtain some $L^2$-results for $\overline\partial$ on forms
that vanish to high order on the singular set of a complex space. As
a  consequence of our main theorem we obtain weighted
$L^2$-solvability results for compactly supported
$\overline\partial$-closed $(p,q)$ forms ($0\le p\le n,\, 1\le q<
n$) on relatively compact subdomains $\Omega$ of the complex space
that satisfy
$H^{n-q}(\Omega,\,\mathcal{S})=0=H^{n-q+1}(\Omega,\,\mathcal{S})$
for every coherent $\mathcal{O}_X$-module $\mathcal{S}$. The latter
result can be used to give an alternate proof of a theorem of Merker
and Porten.
\end{abstract}

\maketitle
\medskip
\noindent \section{Introduction}

\medskip
\noindent This paper is a natural generalization  of results
obtained in an earlier paper of ours \cite{FOV}. There, we addressed
the question of whether one can solve the equation
$\overline\partial u=f$ on $\text{Reg}\,\Omega$, the set of smooth
points of an open relatively compact {\it{Stein}} subset $\Omega$ of
a pure $n$-dimensional reduced {\it{Stein}} space $X$, for a
$\overline\partial$-closed $(p,q)$ form $f$ on $\text{Reg}\,\Omega$
that vanishes to ``high order'' on the singular set of $X$. We
showed that given any $N_0$ non-negative integer there exists a
positive integer $N=N(N_0,\,p,\,q,\,n,\,\Omega)$ such that for any
$\overline\partial$-closed form $f$ on $\text{Reg}\,\Omega,$
vanishing to order $N$ on $\text{Sing}\,X,$ there exists a solution
$u$ to $\overline\partial u=f$ that vanishes to order $N_0$ on
$\text{Sing}\,X$.


\medskip
\noindent In this paper we consider a pure $n$-dimensional countable
at infinity reduced complex space $X$ (we no longer assume that $X$
is Stein). Let $A$ be a nowhere dense, closed analytic subset of
lower dimension containing the singular set $\text{Sing} X$ of $X$
and let $\Omega$ be an open, relatively compact subdomain of $X$
(again we do not assume that $\Omega$ is Stein). We give
$\text{Reg}\,X$ a metric that is compatible with local embeddings
and let $|\;|_x$ and $dV_x$ denote the corresponding pointwise norm
and volume element. Using a partition of unity argument, we can find
a non-negative continuous function $d_A$ on $X$ whose zero set is
$A$ with the property that when $X\supset V^{\text{open}}
\overset{\theta}\to \mathbb{C}^T$ for some positive integer $T,$ is
a local embedding, and $K$ is a compact set in $V$, we have for all
$x\in K$

$$d_A(x)\cong \text{dist}\,(\theta(x),\theta(A)).$$

\noindent Here the right-hand side distance is defined using the
Euclidean metric in $\mathbb{C}^T$ (see section 2 for details on
these constructions).

\medskip
\noindent Let us fix $p\in \mathbb{N}$ with $0\le p\le n$ and let
$\mathcal{H}^{s,\text{loc}}_N(\Omega):=\{f\in
L^{2,\text{loc}}_{p,\,s}(\Omega\setminus A):\,\int_{V\setminus A}
|f|^2 d_A^{-N}dV  \text{\;is\; finite\; for\; all}\;\\
V^{\text{open}}\Subset \Omega \}$. Here $s$ is such that $0\le s\le
n$ and $N$ is a non-negative integer. The main result in this paper
is the following theorem:

\begin{theorem} Let $X$ be a pure $n$-dimensional (countable at infinity)
 reduced complex analytic space. Let $q$ be a positive integer with $1\le q\le n$, and let
 $\Omega$ be a relatively compact  subdomain of $X$ such that $H^q(\Omega,\,\mathcal{S})=0$ for all
 coherent $\mathcal{O}_X$ modules $\mathcal{S}$. For every
non-negative integer $N_0$, there exists a positive integer $N>N_0$
such that: if $f\in \mathcal{H}^{q,\text{loc}}_{N}(\Omega)$ with
$\overline\partial f=0$ on $\Omega\setminus A$ then there exists
$u\in \mathcal{H}^{q-1,\text{loc}}_{N_0}(\Omega)$ with
$\overline\partial u=f$ on $\Omega\setminus A$.
\end{theorem}

\medskip
\noindent Using the above theorem  we obtain weighted
$L^2$-solvability results for compactly supported
$\overline\partial$-closed forms defined on $\Omega\setminus A$.
More precisely we show:

\begin{theorem} Let $X$ be as in Theorem 1.1 and let $q$ be a positive integer
such that $1\le q<n$. Let $\Omega$ be an open relatively compact
subset of $X$ such that
$H^{n-q}(\Omega,\,\mathcal{S})=0=H^{n-q+1}(\Omega,\,\mathcal{S})$
for every coherent $\mathcal{O}_X$-module $\mathcal{S}$.  Let $f$
be a $(p,q)$ form defined on $\Omega\setminus A$ and
$\overline\partial$-closed there, compactly supported in $\Omega$
 and such that $\int_{\Omega\setminus A} |f|^2\, d_A^{N_0} dV<\infty$
for some $N_0\ge 0$. Then there exists a $u\in
L^{2,loc}_{p,\,q-1}(\Omega\setminus A)$ with $\text{supp}_X\,u$
compact in $\Omega$ satisfying $\overline\partial u=f$ on
$\Omega\setminus A$ and such that $\int_{\Omega\setminus A}
|u|^2\, d_A^{N} dV<\infty$, where $N$ is a positive constant that
depends on $N_0$ and $\Omega$.
\end{theorem}

\medskip
\noindent As a corollary of these solvability results we can obtain
an analytic proof of the following theorem of Merker and Porten:

\begin{theorem} (Theorem 2.2 in \cite{MP}) Let $X$ be a connected $(n-1)$-complete normal space
of pure dimension $n\ge 2$. Then, for every domain $D\subset X$ and
every compact set $K\subset D$ with $D\setminus K$ connected,
holomorphic functions on $D\setminus K$ extend holomorphically and
uniquely to $D$.
\end{theorem}

\medskip
\noindent Let us point out here that Merker and Porten prove in
\cite{MP} an extension theorem for meromorphic sections as well. It
is not clear to us at the moment how to employ
$\overline\partial$-techniques to attack the extension problem for
such sections. In November of 2008, Coltoiu \cite{C} and Ruppenthal
\cite{Rup1} independently obtained a $\overline\partial$-theoretic
proof of a Hartogs extension theorem on cohomological
$(n-1)$-complete (resp. $(n-1)$-complete) spaces. The key ingredient
in the proof is a vanishing of the higher direct images of the sheaf
of canonical forms of an appropriate desingularization $\tilde{X}$
of the $(n-1)$-complete complex space $X$. This subtle vanishing
result which was obtained by Takegoshi \cite{T}, easily yields the
vanishing of $H^1_c( \tilde{X}, \mathcal{O})$ which is needed for
the Ehrenpreis method to carry over. Our analytic approach is based
on more general weighted $L^2$-solvability results for
$\overline\partial$-closed, compactly supported $(p,q)$-forms with
$p\neq 0$ which are of independent interest.

\medskip
\noindent The organization of the paper is as follows: In section 2,
we recall some of the key lemmata and theorems from \cite{FOV} and
prove the analogous statements under the more general conditions of
theorem 1.1. In Section 3, we prove Theorem 1.1. Section 4, contains
the weighted $L^2$-solvability results for
$\overline\partial$-closed compactly supported forms. In section 5,
we outline the analytic proof of the generalized Hartogs extension
theorem of Merker-Porten (Theorem 1.3). Section 6 provides a
generalization of Lemma 2.1 in \cite{FOV}. In the short appendix we
prove certain lemmata that are used in the construction of $d_A$.

\section{Preliminaries}

\medskip
\noindent When $X,\Omega$ are Stein ($1$-complete), Theorem 1.1 was
proved in \cite{FOV}. The proof proceeded along the following lines:

\medskip
\noindent {\bf{Step 1: Desingularization.}}  We considered
 $\pi: \tilde{X}\to X$ be a holomorphic surjection with the
following properties:

\smallskip
\noindent i) $\widetilde{X}$ is an $n-$dimensional complex
manifold,\\
ii) $\tilde{A}=\pi^{-1}(A)$ is a hypersurface in $\tilde{X}$ with
only ``normal crossing singularities", i.e. near each $x_0\in \tilde
A$ there are local holomorphic coordinates $(z_1,\dots,z_n)$ in
terms of which $\tilde{A}$ is given by
$h(z)=z_1\cdots z_m=0,$ where $1 \leq m \leq n,$\\
iii) $\pi: \widetilde{X}\setminus \tilde{A} \rightarrow X
\setminus A$ is  a biholomorphism,\\
iv) $\pi$ is proper.

\smallskip
\noindent The existence of such a map follows from  the  facts that
a)\,every reduced, complex space can be desingularized and, b) every
reduced, closed complex subspace of a complex manifold admits an
embedded desingularization (the exact statements and proofs can be
found in \cite{AHV}, \cite{BM}).

\smallskip
\noindent Let $\tilde{\Omega}:=\pi^{-1}(\Omega)$. We gave
$\tilde{X}$ a  real analytic metric $\sigma$ (since by Grauert's
result--Theorem 3 in \cite{Gr2}-- any real analytic manifold can be
properly and analytically embedded in some $\mathbb{R}^T$) and we
considered the corresponding distance function
$d_{\tilde{A}}(x)={\mbox{dist}}(x,\tilde{A}),$ volume element $d
\tilde{V}_{\sigma}$ and norms on $\Lambda^{\cdot} T \tilde{X}$ and
$\Lambda^{\cdot} T^*\tilde{X}.$  Let $J$ denote the ideal sheaf of
$\tilde{A}$ in $\tilde{X}$ and $\Omega^p$  the sheaf of holomorphic
$p$ forms on $\tilde{X}$.  We considered some auxiliary sheaves
(denoted by $\mathcal{L}_{p,q}$)\;on $\tilde{X}$. For every open
subset $U$ of $\tilde{X}$, let \;$\mathcal{L}_{p,q}(U)$ be

\begin{equation}\label{eq:dfn}
\mathcal{L}_{p,q}(U):=\{u\in L^{2,\,{\rm loc}}_{p,q}(U);\;\;
 \overline{\partial}u\in  L^{2,\,{\rm loc}}_{p,q+1}(U)\}
\end{equation}

\noindent  and for each open subset $V\subset U$, let $r^{U}_{V}:
\mathcal{L}_{p,q}(U)\to \mathcal{L}_{p,q}(V)$ be the obvious
restriction maps. Then the map $u\to \overline\partial u$ defines an
$\mathcal{O}_{\tilde{X}}$-homomorphism $\overline\partial:
\mathcal{L}_{p,q}\to \mathcal{L}_{p,q+1}$ and the sequence

$$
0\to \Omega^p\to \mathcal{L}_{p,0}\to \mathcal{L}_{p,1}\to \cdots
\to \mathcal{L}_{p,n}\to 0
$$
\noindent is exact by the local Poincar\'e lemma for
$\overline\partial$. Since each $\mathcal{L}_{p,q}$ is closed under
multiplication by smooth cut-off functions we have a fine resolution
of $\Omega^p$. In the same way, since $J$ is locally generated by
one function, then the sequence

\begin{equation}\label{eq:fres}
0\to J^k \Omega^p\to J^k \mathcal{L}_{p,0}\to \cdots \to J^k
\mathcal{L}_{p,n}\to 0
\end{equation}

\noindent is a fine resolution of $J^k \Omega^p$. Here, $u\in (J^k
\mathcal{L}_{p,q})_{x}$ if it can locally be written as $h^k
\,u_{0}$ where $h$ generates $J_x$  and $u_0\in
(\mathcal{L}_{p,q})_x$. So we can interpret the sheaf cohomology
groups $H^q(\tilde \Omega, (J^k \Omega^p)_{|_{\Omega}})$ as

\begin{eqnarray}\label{eq: acyclicresolution}
H^q(\tilde \Omega, (J^k \Omega^p)_{|_{\Omega}})\cong \frac{
\text{ker}( \overline\partial: J^k
\mathcal{L}_{p,q}(\tilde\Omega)\to J^k
\mathcal{L}_{p,q+1}(\tilde\Omega))}{\text{Im}(\overline\partial: J^k
\mathcal{L}_{p,q-1}(\tilde{\Omega}) \to J^k
\mathcal{L}_{p,q}(\tilde{\Omega}))}.
\end{eqnarray}

\medskip
\noindent{\bf{Step 2: Comparison estimates for pointwise norms of
forms and their pullbacks under the desingularization map.}}  Using
Lojasiewicz inequalities we proved the following pointwise
estimates:

\begin{lemma}(Lemma 3.1 in \cite{FOV}) We have for $x\in \widetilde{\Omega}\setminus \tilde{A},\;
v\in \wedge^{r}T_x (\widetilde{\Omega})$

\begin{eqnarray} \label{eq:pointwise}
c'\; d^{t}_{\tilde{A}}(x)&\le& d_{A}(\pi(x))\le C'\; d_{\tilde{A}}(x),\\
c\; d^{M}_{\tilde{A}}(x) \;|v|_{x,\sigma} &\le&
|\pi_{*}(v)|_{\pi(x)}\le C\; |v|_{x,\sigma}
\end{eqnarray}\\
\noindent  for some  positive constants $c',c,C',C,t,M,$ where
$c,C,M$ may depend on $r$.

\medskip
\noindent For an $r$-form $a$ in $\Omega\setminus A$ set\;
$|\pi^{*}a|_{x,\sigma}:= {\mbox{max}} \{\;|< a_{\pi(x)},
\pi_{*}v>|\;;\; |v|_{x,\,\sigma} \le 1,\; v\in
\wedge^{r}T_{x}(\tilde{\Omega}\setminus \tilde{A})\}$, where by
$<,>$ we denote the pairing of an $r$-form with a corresponding
tangent vector. Using (5) we obtain:

\begin{equation}\label{eq:pointwiseest}
c\; d^{M}_{\tilde{A}}(x) \; |a|_{\pi(x)}\le |\pi^{*}
a|_{x,\sigma}\le C\;|a|_{\pi(x)}
\end{equation}

\noindent on $\tilde{\Omega}$, for some positive constant $M$.
\end{lemma}

\medskip
\noindent {\bf{Step 3: Cohomological vanishing results.}} Inspired
by Grauert's Satz 1, Section 4 in \cite{Gr1}\;(Grauert's result
corresponds to the case where $A$ is a finite set)\; we were led to
the vanishing of a canonical morphism between certain sheaf
cohomology groups on the desingularized space. More precisely we
were able to show the following:

\begin{proposition}(Proposition 1.3 in \cite{FOV})
For $q>0$ and $k \geq 0$ given, there exists a natural number
$\ell,\,\ell\geq k$ such that the map
$$i_*:H^q(\tilde{\Omega}, J^\ell \Omega^p)
 \rightarrow H^q(\tilde{\Omega},  J^k \Omega^p),
$$
\noindent induced by the inclusion $i: J^\ell\,\Omega^p\to
J^k\,\Omega^p$, is the zero map.
\end{proposition}

\noindent The above Proposition played a key role in the proof of
Theorem 1.1, when $X,\Omega$ were Stein. Here is an outline for the
proof of that theorem. Given $N_0$ we chose appropriately $k,\,N$
($k\ge M+t\frac{N_0}{2},\;\;N\ge 2n \ell+M_1$ where $M$ and $t$ are
the exponents that appear in formulas (4) and (6) of lemma 2.1,
$\ell$ is the integer that appears in Proposition 2.2 and $M_1$ is
the exponent that appears on the left-hand side of formula (6) for
the pull-back of the volume form). Then using a change of variables
formula (lemma 4.1 in \cite{FOV}) we had:

\begin{eqnarray}\label{eq:fakepullback}
\int_{\tilde{\Omega}}|\pi^*f|_{\sigma}^2
d_{\tilde{A}}^{-N_1}d\tilde{V}_{\sigma}
 \leq C \int_{\Omega\setminus A}|f|^2d_A^{-N}dV,
\end{eqnarray}

\noindent for a suitable $0<N_1:=N-M_1<N$  and as a consequence we
obtained that $\overline{\partial}\pi^*f=0$ on $\tilde{\Omega}$.

\smallskip
\noindent Using (\ref{eq:fakepullback}) and the way we had chosen
$N, k$ we  showed that $\pi^* f \in J^{\ell}
\mathcal{L}_{p,q}(\tilde \Omega)$. By Proposition 2.2, this implied
that the equation $\overline\partial  v=\pi^*f$ had a solution in
$J^k \mathcal{L}_{p,q-1}(\tilde\Omega)$. Since $|h(x)|\le C
d_{\tilde{A}}(x)$ on compacts in the set where $h$ generates $J$
it  followed that
$$
\int_{\widetilde\Omega'} |v|_{\sigma}^2 d_{\tilde{A}}^{-2k}(x)
d\tilde{V}_{\sigma}<\infty
$$

\noindent where $\widetilde{\Omega}'=\pi^{-1}(\Omega')\;\text
{and}\;\Omega'\subset\subset \Omega$.
\medskip
\noindent Then $\overline\partial((\pi^{-1})^{*} v)=f$ on
$\Omega\setminus A$ and the final step was  to show that

\begin{eqnarray}\label{eq:solutionpseudoestimate}
\int_{\Omega'\setminus A}  |(\pi^{-1})^{*}v|^2 d_{A}^{-N_0} dV \le c
\int_{\widetilde{\Omega'}} |v|_{\sigma}^2 d_{\tilde{A}}^{-2k}
d\tilde{V}_{\sigma}.
\end{eqnarray}

\medskip
\noindent
\subsection{Necessary modifications for the general case.} When
$X,\,\Omega$ are no longer Stein (and hence can not be thought of as
subvarieties embedded in some $\mathbb{C}^T$) we have to work with
locally finite Stein coverings of $X,\Omega$ and carefully modify
the previous approach.

\medskip
\noindent{\bf{Notation:}} Let $a,\;b$ be two non-negative functions
defined on a set $E$. We shall write $a\sim_E b$ iff there exists a
positive constant $C\ge 1$ such that $a\le C\,b$ and $b\le C\,a$ on
$E$.

\medskip
\noindent Clearly $\sim_E$ is an equivalence relation on the set of
non-negative functions defined on $E$. Also, it is not hard to show
that if $\{E_j\}_{j=1}^{M}$ is a finite covering of $E$ and $f,\,g$
are non-negative functions on $E$ then $f\sim_E g$ iff
$f\sim_{E_j}g$ for all $j=1,\cdots,M$.

\medskip
\noindent {\bf{Construction of the distance function $d_A$.}}

\begin{proposition} Let $X$ be a complex space and $A$ be a
subvariety containing $\text{Sing}\,X$. Then, there exists a
non-negative function $d_A$ defined on $X$, with zero set $A$ that
satisfies the following:

\medskip
\noindent If $\phi:U\to C^{T}$ is a local embedding of\;
$U^{\text{open}}\subset X$ and $K,\,K'$ compact subsets of \;$U$
with $K\Subset K'\Subset U$, we have

\begin{equation}\label{eq:dfdist}
d_A(x)\sim_K \text{dist}(\phi(x),\,\phi(K'\cap A))
\end{equation}

\noindent where the right-hand side distance is the Euclidean
distance from the point $\phi(x)$ to the closed set $\phi(K'\cap A)$
in $\mathbb{C}^T$.
\end{proposition}

\medskip
\noindent Before we prove Proposition 2.3 let us make some remarks:
The right-hand side of $(\ref{eq:dfdist})$ can be thought of as the
local definition of $d_A$. There are two things one needs to worry
about: a) How does this local definition depend on the compact
$K'\Supset K$ that contains $K$? b) How does this definition depend
on the embedding? With a little bit of work (see section 7) one can
check that the distances we would get by choosing a different
compact $L$ such that $L\Supset K$ or a different embedding would be
equivalent to each other over $K$.

\medskip
\noindent \noindent{\it{Proof:}} We choose local embeddings $\phi_i:
U'_i\to \mathbb{C}^{T_i}$ where $U'_i$ are open charts in $X$. Let
$U_i$ be open relatively compact subsets of $U'_i$  for
$i=1,2\cdots$ with $U_i\cap A\neq \emptyset$ and such that
$\{U_i\}_{i\ge 1}$ form a locally finite covering of $A$. We choose
an open subset $U_0$ of $X$ with $\overline{U}_0\cap A=\emptyset$
and such that $\{U_i\}_{i\ge 0}$ form a locally finite covering of
$X$. Let $\{\psi_i\}_{i\ge 0}$ be a smooth partition of unity
subordinate to the covering $\{U_i\}_{i\ge 0}$ (i.e. $0\le \psi_i\le
1, \;\;\text{supp}\,\psi_i\subset U_i$ for all $i$ and such that
$\sum_{i=0}^{\infty} \psi_i=1$).

\medskip
\noindent Then we define

\begin{equation}\label{eq:ddist}
d_A(x):=\psi_0(x)+\sum_{i=1}^{\infty} \psi_i(x)
\;\text{dist}(\phi_i(x), \phi_i (\overline{U}_i\cap A)).
\end{equation}

\noindent The i-th right-hand side term in $(\ref{eq:ddist})$ is
extended by $0$ outside $\text{supp}\,\psi_i$. Clearly  $d_A$ is a
non-negative, continuous function on $X$ whose zero set
$d_A^{-1}(0)=A$.

\medskip
\noindent It remains to show that $d_A$ satisfies
$(\ref{eq:dfdist})$. Let $U,\,K,\,K'$ be as in Proposition 2.3. Set
$K_j:=K\cap \text{supp}\,\psi_j$ and $K'_j:=\overline{U}_j\cap K'$.
Since $\{U_j\}_{j\ge 0}$ is a locally finite covering of $X$ there
would exist finitely many $j's$ such that $K\cap
\text{supp}\,\psi_j\neq \emptyset$. If $x\in K\cap
\text{supp}\,\psi_0$, we have

\begin{equation}\label{eq:f1}
0<a_0\le \text{dist}(\phi(x),\,\phi(K'\cap A))\le A_0.
\end{equation}

\medskip
\noindent From the above inequalities one can conclude that when
$x\in K\cap \text{supp}\,\psi_0$ we have that

\begin{equation}\label{eq:f2}
 \psi_0(x)\le 1\le
a_0^{-1}\,\text{dist}(\phi(x),\,\phi(K'\cap A)).
\end{equation}

\medskip
\noindent On the other hand, if $x\in K\cap \text{supp}\,\psi_j$
with $j>0$ we have the following:

\begin{eqnarray*}\label{eq:f3}
\text{dist}(\phi(x),\,\phi(K'\cap A))&\sim_{K_j}&
\text{dist}(\phi(x),\,\phi(K'_j\cap A))\\
&\sim_{K_j}&\text{dist}(\phi_j(x),\,\phi_j(K'_j\cap A))\;\;\;(\text{as\;different\;embeddings\;yield\;equivalent\;norms})\\
&\sim_{K_j}&\text{dist}(\phi_j(x),\,\phi_j(\overline{U}_j\cap A)).\\
\end{eqnarray*}

\medskip
\noindent Let $\{C_j\}_{j>0}$ denote the positive constants that
arise from the fact that

$$\text{dist}(\phi(x),\,\phi(K'\cap A))\sim_{K_j}
\text{dist}(\phi_j(x),\,\phi_j(\overline{U}_j\cap A))$$

\medskip
\noindent i.e. $C_j\ge 1$ and for all $x\in K_j$ we have
$\text{dist}(\phi(x),\,\phi(K'\cap A))\le C_j\,
\text{dist}(\phi_j(x),\,\phi_j(\overline{U_j}\cap A))$ and \\
$\text{dist}(\phi_j(x),\,\phi_j(\overline{U}_j\cap A))\le C_j
\,\text{dist}(\phi(x),\,\phi(K'\cap A))$.

\medskip
\noindent Then, for $x\in K$ we have:

\begin{eqnarray*}
d_A(x)&=&\psi_0(x)+\sum_{i=1}^{\infty} \psi_i(x)
\;\text{dist}(\phi_i(x), \phi_i (\overline{U}_i\cap A))\\
&\le& a_0^{-1}\,\text{dist}(\phi(x),\,\phi(K'\cap
A))+\sum_{\psi_j(x)\neq 0}
C_j\,\psi_j(x)\,\text{dist}(\phi(x),\,\phi(K'\cap A))\\
&\le& C\,\text{dist}(\phi(x),\phi(K'\cap A))
\end{eqnarray*}

\medskip
\noindent where $C:=\text{max}\{a_0^{-1}, A_0, \underset{j>0;
K\cap\text{supp}\,\psi_j\neq \emptyset}\sum C_j\}$. On the other
hand for $x\in K$ we have:

\begin{eqnarray*}
\text{dist}(\phi(x),\,\phi(K'\cap A))&=&\sum
\psi_j(x)\,\text{dist}(\phi(x),\,\phi(K'\cap A))\\
&\le& A_0\,\psi_0(x)+\underset{j>0;\psi_j(x)\neq 0} \sum C_j
\psi_j(x)\,\text{dist}(\phi_j(x),\phi_j(\overline{U}_j\cap A))\\
&\le& C\,d_A(x).
\end{eqnarray*}

\subsection{A Hermitian metric on $\text{Reg}\,X$.} We give
$\text{Reg}\,X$ a hermitian metric $\tilde{\sigma}$ compatible with
local embeddings and let $|\;|_{\tilde{\sigma}}$ denote the
pointwise norms on tangent vectors, cotangent vectors induced by
$\sigma$. Compatible with local embeddings means the following: If
$\phi:U\to \mathbb{C}^T$ is a local embedding of $X$ and $K$ a
compact subset of $U$  we have for all $v\in T_z(\text{Reg}X)$  with
$z\in K\cap\text{Reg}\,X$

$$|\phi_*(v)|_E\sim |v|_{\tilde{\sigma}}.$$

\noindent By $|\;|_E$ we mean the pointwise norm induced by the
Euclidean metric in $\mathbb{C}^T$.

\medskip
\noindent \subsection{Construction of special cut-off functions} In
section 4, we shall need the existence of special cut-off functions.
Using the same notation as above we have:

\begin{proposition} For each $\nu\in \mathbb{N}$ there exists a function $\chi_{\nu}\in
C^{\infty}(X; [0,1])$ such that $\chi_{\nu}(x)=1$ when $d_A(x)\ge
\frac{1}{\nu}$ and $\chi_{\nu}(x)=0$ near $A$. Furthermore for every
$K$ compact subset of $X$ there exists a constant $C_K$ such that

$$
|\overline\partial \chi_{\nu}(x)|_{\tilde{\sigma}}\le C_K\,\nu,
$$

\noindent for all $x\in K\cap \text{Reg}\,X$.
\end{proposition}

\medskip
\noindent{\it{Proof:}} The proof is based on the following well
known principle:

\begin{lemma} For every $d\in \mathbb{N}$, there exists a positive
constant $C_d$ with the following property: For every closed subset
$F$ of $\mathbb{R}^d$ and $\nu\in \mathbb{N}$ there exists
$\chi_{\nu,F}\in C^{\infty}(\mathbb{R}^d; [0,1])$ such that
$\chi_{\nu,F}(x)=1$ when $\text{dist}(x,F)\ge \frac{1}{\nu}$,\;
$\chi_{\nu,F}=0$ near $F$ and $|d\chi_{\nu,F}|_E\le C_d\,\nu$.
\end{lemma}

\medskip
\noindent We consider the family $(\phi_j, U'_j, U_j,\,\psi_j)$ from
the proof of Proposition 2.3 in section 2.1. For $j>0$ taking as
$K:=\text{supp}\psi_j,\;
K':=\overline{U}_j,\;U:=U'_j,\;\phi:=\phi_j$ we obtain from
$(\ref{eq:dfdist})$ the existence of positive integers $m_j$ such
that whenever $\text{dist}(\phi_j(x),\,\phi_j(\overline{U_j}\cap
A))<\frac{\delta}{m_j}\Rightarrow  d_A(x)<\delta$ on
$\text{supp}\,\psi_j$.

\medskip
\noindent For each $\nu,\,j\in \mathbb{N}$ we define (using lemma
2.5)

$$
\chi_{j,\,\nu}:=\chi_{m_j\,\nu,\,\phi_j(\overline{U_j}\cap A)}
$$

\noindent on $\mathbb{C}^{T_j}$ and we set

\begin{equation}\label{eq:dco}
\chi_{\nu}(x):=\psi_0(x)+\sum_{j=1}^{\infty}
\psi_j(x)\,\chi_{j,\,\nu}(\phi_j(x))
\end{equation}

\medskip
\noindent Clearly $\chi_{\nu}\in C^{\infty}(X;[0,1])$ and
$\chi_{\nu}=0$ near $A$. Whenever $d_A(x)\ge \frac{1}{\nu}$, we must
have $\chi_{j,\nu}(x)=1$ when $x\in \text{supp}\,\psi_j$, since
$\text{dist}(\phi_j(x),\,\phi_j(\overline{U_j}\cap A))\ge
\frac{1}{m_j\,\nu}$. Hence
$\chi_{\nu}(x)=\sum_{j=0}^{\infty}\psi_j(x)=1$. Moreover,

$$
\overline\partial \chi_{\nu}=\overline\partial
\psi_o+\sum_{j=1}^{\infty} \left(\overline\partial
\psi_j\wedge\chi_{j,\nu}+\psi_j\,\overline\partial
\chi_{j,\nu}\right).
$$

\medskip
\noindent Given any smooth function $f$ on $X$ we have that
$|\overline\partial f|_{\tilde{\sigma}}$ is bounded on $K\cap
\text{Reg}\,X$ for $K$ compact subset of $X$. Also, given any $K$
compact subset of $X$ there are finitely many $j$'s such that $K\cap
\text{supp}\,\psi_j\neq \emptyset$. Hence for $x\in K\cap
\text{Reg}\,X$ we have

\begin{eqnarray}
|\overline\partial \chi_{\nu}|_{\tilde{\sigma}} &\le&
C'_K+\underset{j;\text{supp}\psi_j\cap K\neq\emptyset}\sum \psi_j
|\phi^*(\overline\partial \chi_{j,\nu})|_{\tilde{\sigma}}\\
&\le& C'_K+\nu C''_K\le (C'_K+C''_K)\,\nu.
\end{eqnarray}

\noindent In the last estimate we used the fact that
$|\overline\partial \chi_{j,\,\nu}|_E\le C_{2T_j}\, m_j\,\nu$ and
that $|\phi_j^*(\theta)|_{\tilde{\sigma}}\le C'_j\,|\theta|_E$ on
$\text{supp}\,\psi_j$.

\subsection{Generalization of results obtained in sections 3 and 4 in \cite{FOV}.} Let $X$ be a pure
$n$-dimensional complex space countable at infinity, let $\Omega$ be
an open relatively compact subdomain of $X$ and let
$\pi:\tilde{X}\to X$ be a desingularization map as in Step 1, of
section 2. Let $d_A$ and $\tilde{\sigma}$ be the distance function
to $A$ and the Hermitian metric on $\text{Reg}\,X$ that were
constructed in sections 2.1-2.2. The pointwise comparison estimates
of forms and their pull-backs under the desingularization map $\pi$
(Lemma 2.1) as well as a change of variables formula that were
obtained in \cite{FOV} carry over to the more general situation we
consider in this paper, using a finite Stein covering of
$\overline{\Omega}$. More precisely we have:

\begin{lemma}
For every $t\in \mathbb{N}^*$, there exists $N\in \mathbb{N}$ such
that if $f\in \mathcal{H}^{s,\,loc}_N(\Omega)$ then $\pi^* f\in
J^t\,\mathcal{L}_{p,s}(\tilde{\Omega})$.
\end{lemma}

\begin{lemma} For every non-negative integer $N_0$, there exists
$t_0\in \mathbb{N}^*$ such that if \;$\tilde{u}\in
J^{t_0}\,\mathcal{L}_{p,s}(\tilde{\Omega})$ then\\
$(\pi^{-1})^*(\tilde{u})\in \mathcal{H}^{s,\,loc}_{N_{0}}(\Omega)$.
\end{lemma}

\medskip
\noindent{\bf{Uniformity of exponents in Lemmata 2.6,\,2.7.}} A
careful inspection of the proofs of Lemmata 3.1 and 4.1 in
\cite{FOV} show that for fixed $p,s,\,N_0,t$ the exponents $N$ and
$t_0$ that appear in Lemmata 2.6,\,2.7 can also be used for forms
$f\in \mathcal{H}^{s,\,loc}_N(W)$ or $\tilde{u}\in
J^{t_0}\,\mathcal{L}_{p,s}(\tilde{W})$ where $W$ is {\it{any}} open
subset of $\Omega$.

\medskip
\noindent
\subsection{Generalization of Proposition 2.2} Let $X$ be a pure $n$-dimensional
complex space countable at infinity and let $\Omega$ be an open
relatively compact subdomain of $X$ and let $\pi:\tilde{X}\to X$ be
a desingularization map as in Step 1, of section 2. Then we have the
following

\begin{proposition} For every $t_0\in \mathbb{N}^*$  there exists a positive integer $t\ge t_0$ such that the map

$$
i_*: R^s\pi_*(J^t\,\Omega^p_{\tilde{X}})\to
R^s\pi_*(J^{t_0}\Omega^p_{\tilde{X}})
$$

\noindent induced by the inclusion homomorphism
$i:J^t\Omega^p_{\tilde{X}}\to J^{t_0}\Omega^p_{\tilde{X}}$, is the
zero homomorphism over $\Omega$, for every $s\in \mathbb{N}^*$.
\end{proposition}

\medskip
\noindent {\it{Proof:}} We can find finitely many Stein subdomains
$\Omega_j\subset\subset \Omega_j^*\subset\subset X$ with
$j=1,\cdots,N_1$ that cover $\Omega$. From the proof of
Proposition 2.2 (Proposition 1.3 in \cite{FOV}) it is clear that
for any $t_0\ge 0$ integer and for any coherent analytic
$\mathcal{O}_{\tilde{X}}$-module $\mathcal{S}$, there exists a
positive integer $t_j\ge t_0$ such that $i_*: H^s
(\tilde{\Omega}_j,\,J^{t_j}\mathcal{S})\to
H^s(\tilde{\Omega}_j,\,J^{t_0}\,\mathcal{S})$ is the zero map for
$s>0$. Take as $\mathcal{S}:=\Omega^p_{\tilde{X}}$. By Satz 5,
Section 2 in \cite{Gr} since each $\Omega_j$ is Stein, we know
that

\begin{equation}\label{eq:Gr}
H^s(\tilde{\Omega}_j,\,J^{t_j}\,\Omega^p_{\tilde{X}})\cong
H^0(\Omega_j,\,R^s\,\pi_* (J^{t_j}\,\Omega^p_{\tilde{X}})).
\end{equation}

\noindent Combining these observations together we conclude that
$H^0(\Omega_j,\,R^s \pi_*(J^{t_j}\,\Omega^p_{\tilde{X}}))\to
H^0(\Omega_j,\,R^s \pi_*(J^{t_0}\Omega^p_{\tilde{X}}))$ is the zero
map. By Cartan's theorem A (Chapter I, vol. III in \cite{GunR}) we
know that for every $x\in \Omega_j$ the germs at $x$ of global
sections $ H^0(\Omega_j,\,R^s\,\pi_*
(J^{t_j}\,\Omega^p_{\tilde{X}}))$ generate the stalk $(R^q
\pi_*(J^{t_j}\Omega^p_{\tilde{X}}))_x$. Hence the induced map $i_*:
R^s\pi_*(J^{t_j}\Omega^p_{\tilde{X}})\to R^s
\pi_*(J^{t_0}\Omega^p_{\tilde{X}})$ is the zero homomorphism over
each $\Omega_j$ since a set of generators of each stalk is mapped to
zero. Let us set $t:=\text{max}(t_1,\,t_2,\cdots,t_{N_1})$. Then
$i_*:R^s\pi_*(J^{t}\Omega^p_{\tilde{X}})\to R^s
\pi_*(J^{t_0}\Omega^p_{\tilde{X}})$ is the zero map over each
$\Omega_j$ and hence over $\Omega$.

\medskip
\noindent If $W$ is a Stein subset of $\Omega$ we see that

$$
H^s(\tilde{W},\,J^{t}\Omega^p_{\tilde{X}})\overset{i_*}\longrightarrow
H^s(\tilde{W},\,J^{t_0}\Omega^p_{\tilde{X}})$$

\noindent is the zero map for all $s>0$. This follows from the
following commutative diagram:

$$
\xymatrix{  &
H^s(\tilde{W},\,J^{t}\Omega^p_{\tilde{X}})\ar[d]_{\cong}
\ar[r]^{i_*} & H^s(\tilde{W},\,J^{t_0}\Omega^p_{\tilde{X}})\ar[d]_{\cong}\\
&H^0(W,\,R^s \pi_*(J^{t}\Omega^p_{\tilde{X}}))\ar[r]^{i_*} &
H^0(W,\,R^s \pi_*(J^{t_0}\Omega^p_{\tilde{X}}))}. $$

\medskip
\noindent The vertical map are isomorphisms by Satz 5, section 2 in
\cite{Gr} and the bottom horizontal map is zero from the result of
the previous paragraph.

\medskip
\noindent{\bf{Uniformity of exponents in Proposition 2.8:}} As
$\Omega_j$ are Stein, by Cartan's theorem A we know that for any
coherent $\mathcal{O}_X$ module $\mathcal{G},$ the space
$H^0(\Omega_j,\mathcal{G})$ generates $\mathcal{G}_x$ for each $x\in
\Omega_j$. If $W_j\subset \Omega_j$ Stein subdomains of $\Omega_j$
for each fixed $s$ we will still have that $i_*:
R^s\,\pi_*(J^{t_j}\mathcal{S})\to R^s\,\pi_*(J^{t_0}\mathcal{S})$ is
the zero map over $W_j$ and $i_*:
H^s(\tilde{W_j},\,J^{t_j}\,\mathcal{S})\to
H^s(\tilde{W_j},\,J^{t_0}\,\mathcal{S})$ is the zero map  for the
same values of $t_j$ (or even $t$) as above.

\medskip

\bigskip \noindent Now using Proposition 2.8, for $q>0$ given and
$t_0$ given  we will inductively define positive integers
$t_1,\,t_2,\cdots,t_q$ with $t_0\le t_1\le \cdots\le t_q$ such
that

\begin{equation}\label{eq:in1}
R^{s} \pi_*(J^{t_{j+1}}\Omega^p_{\tilde{X}})\to
R^{s}\pi_*(J^{t_j}\Omega^p_{\tilde{X}})
\end{equation}

\noindent induced by the inclusion map $i:
J^{t_{j+1}}\Omega^p_{\tilde{X}}\to J^{t_j}\Omega^p_{\tilde{X}}$,
is the zero map over $\Omega$ for $s>0$ and for all $j$ with $0\le
j\le q-1$.

\medskip
\noindent

\section{Proof of Theorem 1.1}

\medskip
\noindent Let $\mathcal{U}=\{U_j\}_{j\in J}$ be a locally finite
Stein covering of $\Omega$ and let
$\tilde{\mathcal{U}}=\{\tilde{U}_j:=\pi^{-1}(U_j)\}_{j\in J}$ be a
locally finite covering of $\tilde{\Omega}=\pi^{-1}(\Omega)$. Let
$r,s,t$ be non-negative integers. We introduce the spaces of
alternating $r$ co-chains

$$C^{r,s}_t:=C^r(\tilde{\mathcal{U}}, J^t\mathcal{L}_{p,s}).$$

\medskip
\noindent Here $\mathcal{L}_{p,s}$ are the sheaves defined in
$(\ref{eq:dfn})$. Let $\delta: C^{r,s}_t\to C^{r+1,s}_t$ be the
co-boundary map and $\overline\partial: C^{r,s}_t\to C^{r,s+1}_t$
the $\overline\partial$ operator applied componentwise, i.e. given
any $c=\{c_{j_0j_1\cdots j_{r}}\}$ we define $\overline\partial c$
to be the $r$ co-chain such that $(\overline\partial
c)_{j_0j_1\cdots j_{r}}=\overline\partial c_{j_0j_1\cdots j_{r}}$.
Clearly $\delta\circ\overline\partial=\overline\partial\circ
\delta$. Set $U_{j_0j_1\cdots j_r}:=U_{j_0}\cap U_{j_1}\cap\cdots
\cap U_{j_r}$ (as $U_j$ are Stein, $U_{j_0j_1\cdots j_r}$ are
Stein).

\medskip
\noindent Let $N_0,\,q$ be as in Theorem 1.1. Using Lemma 2.7 we
can find an integer $t_0$ such that if \;$\tilde{u}\in
J^{t_0}\,\mathcal{L}_{p,q-1}(\tilde{\Omega})$ then
$(\pi^{-1})^*(\tilde{u})\in\mathcal{H}^{q-1,\,loc}_{N_{0}}(\Omega)$.
With $t_0$ chosen as before let us choose $\{t_j\}_{1\le j\le q}$,
positive integers with $t_0\le t_1\le \cdots \le t_q$ in such a
way such that the maps determined by $(\ref{eq:in1})$  are the
zero maps over $\Omega$ for all $0\le j<q$.

\medskip
\noindent \begin{claim} Given any $c\in C^{r,s}_{t_{j+1}}$ with
$s>0$ and $\overline\partial c=0$ we can find a solution to
$\overline\partial d=c$ with $d\in C^{r,s-1}_{t_j}$ for $0\le
j<q$.
\end{claim}

\medskip
\noindent{\it{Proof of Claim:}} The Claim follows easily from Step
1 in section 2.1 taking into account the facts that
$U_{j_0j_1\cdots j_r}$ are Stein, $(\ref{eq: acyclicresolution})$
holds and that the maps in $(\ref{eq:in1})$ are zero over $\Omega$
for all $j$ with $0\le j\le q-1$. $\phantom{sd}\square$

\medskip
\noindent Now, with $q$ and $t_q$ chosen as above, Lemma 2.6
determines an integer $N$ such that if $f\in
\mathcal{H}^{q,\,loc}_ {N}(\Omega)$ with $\overline\partial f=0$
on $\Omega\setminus A$ then $\pi^*f\in
J^{t_q}\,\mathcal{L}_{p,q}(\tilde{\Omega})$ and $\overline\partial
\pi^*f=0$ on $\tilde{\Omega}$.

\medskip
\noindent The proof of Theorem 1.1 will proceed as follows:

\medskip
\noindent {\bf{Part A:}} To $\pi^*f$ we will associate a $q$
co-cycle $c^q\in
Z^q(\tilde{\mathcal{U}},\,J^{t_0}\Omega^p_{\tilde{X}})$ (in a
similar way as when we prove the quasi-isomorphism between Dolbeault
cohomology and \v Cech-cohomology on a manifold).

\medskip
\noindent Using $\pi^*f$ we define a $0$ co-cycle $c^0:=\{c^0_j\}$
such that $c^0\in C^{0,q}_{t_q}$ by setting $c^0_j:=\pi^*
f_{\upharpoonright (\tilde{U_j}\setminus \tilde{A})}$ for all $j\in
\mathbb{N}$. Clearly $\delta c^0=0$ and $\overline\partial c^0=0$.

\medskip
\noindent \medskip \noindent Using Claim 3.1, we can inductively
define $r$ co-chains $c^r\in C^{r,q-r}_{t_{q-r}}$ for $1\le r\le
q$ and $d^r\in C^{r,q-r-1}_{t_{q-r-1}}$ for $0\le r<q$ such that

\medskip
\noindent(i) $\overline\partial d^r=c^r$ for all $r$ with $0\le
r\le q-1$,

\medskip
\noindent (ii) $c^{r+1}=\delta d^r$ for all $r$ with $0\le r\le
q-1$,

\medskip
\noindent (iii) $\delta c^r=0\;\;\text{and}\;\;\overline\partial
c^r=0$ for all $r$ with $0\le r\le q$.

\medskip
\noindent Suppose that $c^0,\,c^1,\cdots, c^r$ and
$d^0,d^1,\cdots, d^r$ with $0\le r\le q-1$ satisfying (i)-(iii)
are given. Set $c^{r+1}:=\delta d^r$. Clearly $\delta c^{r+1}=0$
and $\overline\partial c^{r+1}=\overline\partial \delta d^r=\delta
\overline\partial d^r=\delta c^r=0$. If $r+1<q$ then by Claim 3.1,
we can find $d^{r+1}$ such that $\overline\partial
d^{r+1}=c^{r+1}$. Then $c^{r+1},\,d^{r+1}$ satisfy (i)-(iii) and
the procedure is completed. If $r+1=q$ then the previous procedure
yielded an element $c^q=\delta d^{q-1}\in C^{q,\,0}_{t_0}$ such
that $\overline\partial c^q=0$. But the latter implies that $c^q$
is a co-cycle in
$C^q(\tilde{\mathcal{U}},\,J^{t_0}\,\Omega^p_{\tilde{X}})$.

\medskip
\noindent {\bf{Part B:}} Taking into account that for any coherent
analytic sheaf $\mathcal{S}$ on $\tilde{X}$ and for any $U$ open in
$X$ we have $\pi_* \mathcal{S}(U)=\mathcal{S}(\tilde{U})$, we can
see that $c^q$ will define a $q$ co-cycle $\overline{c}^q$ in
$C^q(\mathcal{U},\,\pi_* J^{t_0}\Omega^p_{\tilde{X}})$. Since
$\mathcal{U}$ is a Stein covering of $\Omega$ we have that
$H^q(\mathcal{U},\,\pi_*(J^{t_0}\Omega^p_{\tilde{X}}))\cong
H^q(\Omega,\pi_*(J^{t_0}\Omega^p_{\tilde{X}}))$. But the latter
group vanishes since $\pi_*(J^{t_0}\Omega^p_{\tilde{X}})$ is a
coherent sheaf on $X$ and by assumption
$H^q(\Omega,\,\mathcal{S})=0$ for all coherent
$\mathcal{O}_X$-modules $\mathcal{S}$. Hence there exists a
$\tilde{h}^{q-1}\in C^{q-1}(\mathcal{U},
\pi_*(J^{t_0}\Omega^p_{\tilde{X}}))$ such that $\delta
\overline{h}^{q-1}=\overline{c}^q$ and thus there exists an
$h^{q-1}\in
C^{q-1}(\tilde{\mathcal{U}},J^{t_0}\Omega^p_{\tilde{X}})$ such that
$\delta h^{q-1}=c^q$. We set $d'^{q-1}:=d^{q-1}-h^{q-1}$. Clearly
$\delta d'^{q-1}=0$ and $\overline\partial d'^{q-1}=c^{q-1}$.

\medskip
\noindent If $q=1$ then we are done. If $q>1$ then  by downward
induction on $r$, we shall modify $d^r$ to $d'^r\in
C^{r,\,q-r-1}_{t_0}$ satisfying for all $r$ with $0\le r\le q-1$

\begin{equation}\label{eq:cc}
\delta d'^r=0\;\;\text{and}\;\;\;\;\overline\partial d'^r=c^r.
\end{equation}

\medskip
\noindent We start with $d'^{q-1}$ as above. Suppose $d'^j$ are
given for $j=q-1,\,q-2,\cdots r$ (with $r>0$)  satisfying
$(\ref{eq:cc})$. As $\delta d'^r=0$ we will construct a
$(q-1)$-co-chain  $a^{r-1}\in C^{r-1,\,q-r-1}_{t_0}$ that satisfies
$\delta a^{r-1}=d'^r$. Let $\{\phi_j\}$ be a smooth partition of
unity subordinate to the cover $\tilde{\mathcal{U}}$ of
$\tilde{\Omega}$. We define an $(r-1)$ co-chain $a^{r-1}\in
C^{r-1,\,q-r-1}_{t_0}$ by setting:

\begin{equation}\label{eq:da}
(a^{r-1})_{j_0j_1\cdots j_{r-1}}:=\sum \phi_j {d'}^r_{jj_0\cdots
j_{r-1}}
\end{equation}

\medskip
\noindent It is a standard fact that $\delta a^{r-1}=d'^r$. Taking
$\overline\partial$ on both sides of the above equation we obtain
$\overline\partial \delta a^{r-1}=\overline\partial d'^r=c^r=\delta
d^{r-1}$. Hence we have $\delta (d^{r-1}-\overline\partial
a^{r-1})=0$. Let us set $d'^{r-1}:=d^{r-1}-\overline\partial
a^{r-1}$. Then $\overline\partial d'^{r-1}=c^{r-1}$ and $\delta
d'^{r-1}=0$. It remains to show that
$d'^{r-1}=d^{r-1}-\overline\partial a^{r-1}\in C^{r-1,q-r}_{t_0}$
for the induction step to be completed. From $(\ref{eq:da})$ we
obtain:

$$
\overline\partial a^{r-1}_{j_0\cdots j_{r-1}}=\sum_j
\overline\partial \phi_j\wedge d'^r_{jj_0\cdots j_{r-1}}+ \sum
\phi_j c^r_{jj_0\cdots j_{r-1}}\in J^{t_0}\mathcal{L}_{p,q-r},
$$

\noindent as $|\overline\partial \phi_j|_{\tilde{\sigma}}$ is
bounded for each $j$ and $\{\text{supp}\,\phi_j\}_j$ is locally
finite. Hence $d'^r\in C^{r-1, q-r}_{t_0}$.

\medskip
\noindent Hence we can find $d'^0\in C^{0,q-1}_{t_0}$ with $\delta
d'^0=0$ and $\overline\partial d'^0=c^0$. But $\delta d'^0=0$
implies that ${d'}^0_j={d'}^0_k$ on $\tilde{U}_j\cap \tilde{U}_k$
(when non-empty). Hence $d'^0\in C^{0,q-1}_{t_0}$ defines a form
$\tilde{u}\in J^{t_0}\mathcal{L}_{p,q-1}(\tilde{\Omega})$ such that
$\overline\partial \tilde{u}=\pi^* f$ on $\tilde{\Omega}\setminus
\tilde{A}$. Then $u:=({\pi^{-1}})^* \tilde{u}\in
\mathcal{H}^{q-1,\,loc}_{N_0}$ (by Lemma 2.7) and $\overline\partial
u=f$ on $\Omega\setminus A$.

\medskip
\noindent

\medskip
\noindent{\bf{Remark:}} If $W\Subset \Omega$ open subdomain of
$\Omega$ such that $H^q(W,\,\mathcal{S})=0$ for all coherent
$O_X$-modules $\mathcal{S}$, then for any $N_0$ non-negative integer
and any $f\in \mathcal{H}^{q,\text{loc}}_{N}(W)$ with
$\overline\partial f=0$ on $W\setminus A$  there exists $u\in
\mathcal{H}^{q-1,\text{loc}}_{N_0}(W)$ with $\overline\partial u=f$
on $W\setminus A$. Here $N$ is the same integer as the one that
appears in Theorem 1.1. This follows from the uniformity of
exponents in Lemmata 2.6, 2.7 and Proposition 2.8, along with the
fact that the proof of Claim 3.1 carries over verbatim for such a
$W$.

\bigskip
\section{Weighted $L^2$-solvability results for
$\overline\partial$-closed compactly supported forms}

\medskip
\noindent

\noindent In \cite{OV} (Section 5) we developed an  analytic
approach to obtaining Hartogs extension theorems on normal Stein
spaces (of pure dimension $n\ge 2$) with arbitrary singularities.
One of the key elements in that approach was obtaining weighted
$L^2$-solvability results for $\overline\partial$-closed compactly
supported forms. More precisely we proved the following (under the
assumption that $X,\Omega$ are $1$-complete spaces and for
$A:=\text{Sing}\,X$):

\begin{theorem} (Theorem 5.3 in \cite{OV}) Let $f$ be a
$(p,q)$ form defined on $\text{Reg}\,\Omega$ and
$\overline\partial$-closed there with $0<q<n$, compactly supported
in $\Omega$ and such that $\int_{\text{Reg}\,\Omega} |f|^2\,
d_A^{N_0} dV<\infty$ for some $N_0\ge 0$. Then there exists a
solution $u$ to $\overline\partial u=f$ on $\text{Reg}\,\Omega$
satisfying $\text{supp}_X\, u\Subset \Omega$ and such that

$$\int_{\text{Reg}\,\Omega} |u|^2\, d_A^{N} dV\le
C\,\int_{\text{Reg}\,\Omega} |f|^2\, d_A^{N_0} dV
$$

\noindent $N$ depends on $N_0$ and $\Omega$ and $C$ is a positive
constant that depends on $N_0,\,N,\,\Omega$\;
and\;\;$\text{supp}\,f$.

\end{theorem}

\medskip
\noindent Let $X$ now be a pure $n$-dimensional complex space
reduced and countable at infinity, $q$ be a positive integer with
$q< n$ and $\Omega\Subset X$ open such that
$H^{n-q}(\Omega,\,\mathcal{S})=0=H^{n-q+1}(\Omega,\,\mathcal{S})$
for every coherent $\mathcal{O}_X$ module $\mathcal{S}$ and let
$A$ be a nowhere dense, lower dimensional complex analytic subset
of $X$ containing $\text{Sing}\,X$. Let $\mathcal{U}:=\{U_i\}$ be
a locally finite Stein cover of $\Omega$, consisting of open
relatively compact Stein subsets of $\Omega$ such that each $U_i$
is relatively compact in $U'_i$ (the local coordinate charts
$(\phi_i,\,U'_i)$ that were described in section 2.2). Let
$\{\psi_i\}$ be a smooth partition of unity subordinate to this
cover. Using local embeddings we define the distance $d_A$ as in
section 2.1. We give $\text{Reg}\,X$ the metric $\tilde{\sigma}$
which was defined in section 2.2. This metric is compatible with
local embeddings and for simplicity let $|\;|_x$ and $dV_x$ denote
the corresponding pointwise norm and volume element with respect
to this metric. Then we can prove the following generalization of
Theorem 4.1:

\begin{theorem} Let $X,\Omega,\,A,\,q$ be as above. Let $f$
be a $(p,q)$ form defined on $\Omega\setminus A$ and
$\overline\partial$-closed there, compactly supported in $\Omega$
 and such that $\int_{\Omega\setminus A} |f|^2\, d_A^{N_0} dV<\infty$
for some $N_0\ge 0$. Then there exists a $u\in
L^{2,loc}_{p,\,q-1}(\Omega \setminus A)$ with $\text{supp}_X\,u$
compact in $\Omega$ satisfying $\overline\partial u=f$ on
$\Omega\setminus A$ and such that $\int_{\Omega\setminus A}
|u|^2\, d_A^{N} dV<\infty$, where $N$ is a positive constant that
depends on $N_0$ and $\Omega$.
\end{theorem}

\medskip
\noindent{\it{Proof:}} The proof follows along the same lines as
the proof of Theorem 4.1 (Theorem 5.3 in \cite{OV}). Let us use
the symbol $\mathcal{H}^{s,\,loc}_N(\Omega):=\{f\in
L^{2,\,loc}_{n-p,\,s}(\Omega\setminus A);\;\;\int_{V\setminus A}
|f|^2 \,d_A^{-N}\,dV<\infty\;\;\text{for\;\,all}\;V\Subset
\Omega\}$. As the domain $\Omega$ satisfies
$H^{n-q}(\Omega,\,\mathcal{S})=0=H^{n-q+1}(\Omega,\,\mathcal{S})$
for all coherent $\mathcal{O}_X$-modules $\mathcal{S}$, we can
apply Theorem 1.1 twice for forms of different bidegree. More
precisely for $N_0+2$ there exists $N_1 (>>N_0+2)$ such that if
$F\in \mathcal{H}^{n-q,\, loc}_{N_1}(\Omega)$, $\overline\partial
F=0$ on $\Omega\setminus A$ then, there exists $a\in
\mathcal{H}^{n-q-1,\,loc}_{N_0+2}(\Omega)$ satisfying
$\overline\partial a=F$ on $\Omega\setminus A$. Similarly, for the
above $N_1$ there exists an $N (>> N_1)$ such that if $G\in
\mathcal{H}^{n-q+1,\, loc}_{N}(\Omega),\;\overline\partial$-closed
on $\Omega\setminus A$ then, there exists a \; $b\in \mathcal{H}^{
n-q,\,loc}_{N}(\Omega)$ satisfying $\overline\partial b=G$ on
$\Omega\setminus A$.

\medskip
\noindent Let $f,\,N_0$ be as in Theorem 4.2 and $N,\,N_1$ be chosen
as above.  Consider the following map:

$$ L_f: \mathcal{H}^{n-q+1,\,loc}_N (\Omega)\cap \text{kern}(\overline\partial)\to
\mathbb{C}$$

\noindent defined by
$$L_f(w):=(-1)^{p+q+1}\;\int_{\Omega\setminus A} v\wedge f$$

\noindent where $v\in \mathcal{H}^{n-q,\,loc}_{N_1}(\Omega)$ is a
solution to $\overline\partial v=w$ on $\Omega\setminus A$ (such a
solution always exist by Theorem 1.1).

\medskip
\noindent First of all we need to show that $L_f$ is well-defined,
i.e. independent of the choice of the solution $v\in
\mathcal{H}^{n-q,\,loc}_{N_1}$ to the equation $\overline\partial
v=w$ on $\Omega\setminus A$. It suffices to show that $\int_{
\Omega\setminus A} v\wedge f=0$ when  $v\in
\mathcal{H}^{n-q,\,loc}_{N_1}(\Omega)$ and $\overline\partial v=0$
on $\Omega\setminus A$. According to what was discussed above
there exists an $a\in \mathcal{H}^{n-q-1,\,loc}_{N_0+2}(\Omega)$
satisfying $\overline\partial a=v$ on $\Omega\setminus A$.

\medskip
\noindent Let $\chi_{\nu}\in C^{\infty}(\text{Reg}\,X)$ be the
smooth cut-off functions that were constructed in section 2.3.
Recall that $0\le \chi_{\nu}\le 1,$\; $\chi_{\nu}(z)=1$ when
$d_A(z)>\frac{1}{\nu},\;\;\;\chi_{\nu}(z)=0$ if $d_A(z)\le
\frac{1}{2\,\nu}$ and such that for any compact $K\subset X$ we
have $|\overline\partial \chi_{\nu}|\le {C_K}{\nu}$ on
$\text{Reg}\,X\cap K$,  for some positive constant $C_K$
independent of $\nu$. Then $L_f(w)=\underset{\nu\to \infty}\lim
\int_{\Omega\setminus A}\chi_{\nu}\, v\wedge f=\underset{\nu\to
\infty}\lim \int_{\Omega\setminus A} \chi_{\nu}\,\overline\partial
a\,\wedge f=-\underset{\nu\to \infty}\lim \int_{\Omega\setminus A}
\overline\partial \chi_{\nu}\wedge a\wedge f$. The last equality
follows from integration by parts arguments; the form
$\chi_{\nu}\,f$ has compact support on $\Omega\setminus A$ but is
not really smooth, so a smoothing argument is needed in order to
apply the standard Stokes' theorem. Let $I_{\nu}:=
\int_{\Omega\setminus A} \overline\partial \chi_{\nu}\wedge
a\wedge f$. By Cauchy-Schwarz inequality we have

\begin{equation*}
|I_{\nu}|^2\le \left(\int_{\text{supp}\,f\cap
\text{supp}\,(\overline\partial\chi_{\nu})}\;
|d_A\,\overline\partial \chi_{\nu}|^2\;|f|^2\;
d_A^{N_0}\;dV\right)\;\left(\int_{\{{\chi_{\nu}}<1\}\cap
\text{supp}\,f}\; |a|^2\;d_A^{-(N_0+2)}\;dV\right)=A\cdot B
\end{equation*}

\smallskip
\noindent But $B$ is finite and  $A\to 0$ as $\nu\to \infty$ since
$|d_A\,\overline\partial \chi_{\delta}|\le C$ and
$\int_{\Omega\setminus A}\;|f|^2\, d_A^{N_0}\;dV<\infty$. Hence
$\underset{\nu\to \infty}\lim I_{\nu}=0$ and thus $L_f$ is
well-defined.

\medskip
\noindent Clearly $L_f$ is a linear map. By Cauchy-Schwarz we have

\begin{equation}\label{eq:CS1}
|L_f(w)|\le C_0\,\left(\int_{\text{supp}\,f} |v|^2\,
d_A^{-N_1}\;dV\right)^{\frac{1}{2}}\;\left(\int_{\Omega\setminus
A}\,|f|^2 \;d_A^{N_0}\;dV\right)^{\frac{1}{2}}
\end{equation}

\noindent To obtain this inequality we used the fact that $N_1>>N_0$
and hence $d_A^{N_1-N_0}\le C_0$ on $\overline{\Omega}$. We want to
show that $L_f$ factors into a bounded linear functional on a
subspace $\mathcal{A}$ of a Hilbert space. Recall that  the
following lemma was proven in \cite{FOV} using the open mapping
theorem for Fr\'echet spaces.

\begin{lemma}(Lemma 4.2 in \cite{FOV}) Let $M$ be a complex manifold and let $E$ and $F$ be Fr\'echet
spaces of differential forms (or currents) of type $(p,q-1),
\;(p,q)$, whose topologies are finer (possibly strictly finer) than
the weak topology of currents. Assume that for every $f\in F$, the
equation $\overline\partial u=f$ has a solution $u\in E$. Then, for
every continuous seminorm $p$ on $E$, there is a continuous seminorm
$q$ on $F$ such that the equation $\overline\partial u=f$ has a
solution with $p(u)\le q(f)$ for every $f\in F,\;q(f)>0$.
\end{lemma}

\medskip
\noindent Let $p(v):=\left(\int_{\text{supp}\,f} |v|^2\,
d_A^{-N_1}\;dV\right)^{\frac{1}{2}}$. Using the lemma in our
situation, given the seminorm $p$ there exists an open set
$W\Subset \Omega$ (that depends on the $\text{supp}\,f$) such that
for all $w\in \mathcal{H}^{n-q+1,\,loc}_N(\Omega)\cap
\text{kern}(\overline\partial)$ with $q(w)=\left(\int_W
|w|^2\,d_A^{-N}\,dV\right)^{\frac{1}{2}}>0$ there exists a
solution $v$ to $\overline\partial v=w$ on $\Omega\setminus A$ and
a positive constant $C$ satisfying

\begin{equation}\label{eq:CS2}
\left(\int_{\text{supp}\,f} |v|^2\,
d_A^{-N_1}\;dV\right)^{\frac{1}{2}}\le C \,q(w).
\end{equation}

\medskip
\noindent If $q(w)=0$ then the same argument will imply that for
every $\epsilon>0$ there exists a solution $v_{\epsilon}$ to
$\overline\partial v_{\epsilon}=w$ on $\Omega\setminus A$ with
$p(v_{\epsilon})<\epsilon$. Then for such a $w$ we have:
$|L_f(w)|\le C_0 \epsilon \int_{\Omega\setminus A}
|f|^2\,d_A^{N_0}\,dV$. Here we used the fact that $L_f(w)$ is
well-defined independent of the choice of solution $v\in
\mathcal{H}^{n-q,\,loc}_{N_1}(\Omega)$. Taking the limit as
$\epsilon\to 0$ we obtain

\begin{equation}\label{eq:CS3}
|L_f(w)|=0\le C_0\,q(w)\,\left(\int_{\Omega\setminus A}\,|f|^2
\;d_A^{N_0}\;dV\right)^{\frac{1}{2}}.
\end{equation}

\medskip
\noindent Combining
$(\ref{eq:CS1}),\,(\ref{eq:CS2}),\,(\ref{eq:CS3})$ we obtain for all
$w\in \mathcal{H}^{n-q+1,\,loc}_N\cap
\text{kern}(\overline\partial)$

\begin{equation}\label{eq:cont}
|L_f(w)|\le C_0\,C\,q(w)\,\left(\int_{\Omega\setminus A}\,|f|^2
\;d_A^{N_0}\;dV\right)^{\frac{1}{2}}.
\end{equation}

\medskip
\noindent From $(\ref{eq:cont})$ we see that  $L_f(w)$ depends
only on $w_{\upharpoonright W}$. Indeed, let $w,\,w'\in
\mathcal{H}^{n-q+1,\,loc}_N\cap \text{kern}(\overline\partial)$
such that $w_{\upharpoonright W}=w'_{\upharpoonright W}$. Then
$L_f(w)=L_f(w-w'+w')=L_f(w-w')+L_f(w')$. From $(\ref{eq:cont})$ we
obtain that $|L_f(w-w')|\le C_0\,C\,
q(w-w')\,\left(\int_{\Omega\setminus A}\,|f|^2
\;d_A^{N_0}\;dV\right)^{\frac{1}{2}}$. But $q(w-w')=0$ as $w-w'=0$
on $W$. Hence $L_f(w-w')=0$ and thus $L_f(w)=L_f(w')$. Hence $L_f$
factors to a well-defined bounded linear functional on
$\mathcal{A}:=\{w_{\upharpoonright W};\;\;w\in
\mathcal{H}^{n-q+1,\,loc}_N(\Omega)\cap
\text{kern}(\overline\partial)\}\subset \mathcal{H}^{n-q+1}_N(W)$.
Here $\mathcal{H}^{n-q+1}_N(W):=\{w\in
L^{2,\,loc}_{n-p,\,n-q+1}(W\setminus A): \;=\;\int_{W\setminus A}
|F|^2\;d_A^{-N}\,dV<\infty\}$.

\medskip
\noindent  We make a norm-preserving extension of the above
functional $L_f$ to $\mathcal{H}^{n-q+1}_N(W)$. Let us call
$\tilde{L}_f$ the extended functional. By Riesz representation
theorem there exists a $u'\in \mathcal{H}^{n-q+1}_N(W)$ such that
for all $w\in \mathcal{H}^{n-q+1}_N(W)$ we have

\begin{equation}\label{eq:Dfn}
\tilde{L}_f(w)=\int_{W\setminus A} <w,\,u'>\;d_A^{-N}\;dV.
\end{equation}

\noindent Set $u:=d_A^{-N}\;\overline{*}\;u'$ on $W\setminus A$
(here $*$ is the Hodge-star operator) and extend by zero outside
$\overline{W}$. We claim that $u$ is the desired solution of
Theorem 5.3. Certainly $\text{supp}\,u\Subset \Omega$ and
$\int_{\Omega\setminus A} |u|^2\,d_A^N\;dV=\int_{W\setminus A}
|u'|^2\; d_A^{-N}\, dV<\infty$, since $u'\in
\mathcal{H}^{n-q+1}_N(W)$. We can control the weighted $L^2$ norm
of $u'$ in terms of the weighted $L^2$-norm of $f$ taking into
consideration the following:

\begin{equation}\label{eq:L2normu'}
\left(\int_{W\setminus A}
|u'|^2\,d_A^{-N}\,dV\right)^{\frac{1}{2}}= \|\tilde{L}_f\|=
\|{L_{f}}_{\upharpoonright
\mathcal{A}}\|=\phantom{sdsdfsdsdsdfhjjlklliuiu}
\end{equation}

\begin{equation*}
=\text{sup}\{ |L_f(w)|:\; w\in \mathcal{A} \;\text{with}\; q(w)\le
1 \}\le C_0\,C\,\left(\int_{\Omega\setminus A}\,|f|^2
\;d_A^{N_0}\;dV\right)^{\frac{1}{2}}.
\end{equation*}

\medskip
\noindent It remains to show that $\overline\partial u=f$ on
$\Omega\setminus A$. Let $\phi\in C^{\infty}_{0,\,(n-p,\,
n-q)}(\Omega\setminus A)$ be a smooth compactly supported form of
bidegree $(n-p,\,n-q)$ on $\Omega\setminus A$. We need to show
that

\begin{equation}\label{eq:Lsolv}
\int_{\Omega\setminus A} \overline\partial \phi\wedge u
=(-1)^{p+q+1}\, \int_{\Omega\setminus A}  \phi \wedge f.
\end{equation}

\medskip
\noindent But $\phi \in \mathcal{H}^{n-q,\,loc}_{N_1}(\Omega)$,
$\overline\partial\phi \in \mathcal{H}^{n-q+1,\,loc}_N(\Omega)$ and
$\overline\partial \phi_{\upharpoonright W}\in \mathcal{A}$.
Therefore from the definition of $L_f$ we have that

$$\tilde{L}_f(\overline\partial \phi_{\upharpoonright
W})=L_f(\overline\partial \phi_{\upharpoonright
W})=(-1)^{p+q+1}\,\int_{\Omega\setminus A} \phi\wedge f$$

\medskip
\noindent On the other hand from $(\ref{eq:Dfn})$ we have that

$$
\tilde{L}_f(\overline\partial \phi_{\upharpoonright
W})=\int_{W\setminus A} \overline\partial \phi \wedge
u=\int_{\Omega\setminus A} \overline\partial \phi \wedge u.
$$

\medskip
\noindent Putting the last two equalities together we obtain
$(\ref{eq:Lsolv})$.

\medskip
\noindent By the definition of $u$ and $(\ref{eq:L2normu'})$ we can
obtain the following estimate for the solution $u$:

$$
\int_{\Omega\setminus A} |u|^2\,d_A^{N}\,dV\le
\tilde{C}\,\int_{\Omega\setminus A} |f|^2\,d_A^{N_0}\,dV
$$

\noindent where $\tilde{C}$ is a positive constant that depends on
$N,\,N_0,\,\Omega,\,\text{supp}\,f$.

\medskip
\noindent

\medskip
\noindent{\bf{Remark:}} We can be more precise about the
dependence of $\text{supp}\,u$ and the constant $\tilde{C}$ (that
appears in the last inequality) on  $\text{supp}\,f$. Let
$X,\Omega$ be as in theorem 4.2 and let $N_0$ be a non-negative
integer. There exists a positive integer $N$ that depends on $N_0$
and $\Omega$ such that the following is true: For every compact
$K\subset \Omega$, there exists a compact $K'\subset \Omega$ and a
positive constant $C$ that depends on $K,\,N,\,N_0,\,\Omega$ such
that for every $(p,q)$ form $f$ with $\text{supp}\,f\subset K$ and
$\int_{K\setminus A}|f|^2\,d_A^{N_0}\,dV<\infty$,
$\overline\partial$-closed on $\Omega\setminus A$, there exists a
solution $u$ to $\overline\partial u=f$ on $\Omega\setminus A$
with $\text{supp}\,u\subset K'$ satisfying

$$
\int_{K'\setminus A} |u|^2\,d_A^{N}\,dV\le
\tilde{C}\,\int_{K\setminus A} |f|^2\,d_A^{N_0}\,dV.
$$

\medskip
\noindent The proof follows along the same lines as the proof of
Theorem 4.2, by taking as $p(v):=\left(\int_{K\setminus A}
|v|^2\,d_A^{-N_1}\,dV\right)^{\frac{1}{2}}$. Then there exist an
open, relatively compact subset $W$ of $\Omega$ (that depends on
$K$) a seminorm $q(w):=\int_{W\setminus A} |w|^2\,d_A^{-N}\,dV$ and
a positive constant $C'$ (that depends on $K$) such that the
equation $\overline\partial v=w$ has a solution $v$ satisfying
$p(v)\le C'\,q(w)$. The rest of the proof of Theorem 4.2 carries
over. The compact $K'$ is chosen to be $K':=\text{closure}(W)$ and
the constant $\tilde{C}=C'\,C_0$ where $C_0$ depends on
$N,\,N_0,\,\Omega$.

\section{An analytic proof of Theorem 1.3}

\medskip
\noindent In \cite{OV} the following proposition was proved:

\begin{proposition}(Proposition 5.1 in \cite{OV}) Let $X$ be a connected, non-compact normal
complex space and $K$ a compact subset of $X$. Let us assume that
$K$ has an open neighborhood $\Omega$ in $X$ with the following
property (P)

\medskip
\noindent For every ``nice'' $\overline\partial$-closed $(0,1)$-form
$f$ on $\text{Reg}\,\Omega$ with
$\text{supp}_X\,f:=\text{closure}_X\,(\text{supp}\,f)$ compact in
$\Omega$, the equation $\overline\partial u=f$ has a solution on
$\text{Reg}\,\Omega$ with
$\text{supp}_X\,u:=\text{closure}_X\,(\text{supp}\,u)$ compact in
$\Omega$.

\medskip
\noindent Then, for every open  neighborhood $D$ of $K$ with
$D\setminus K$ connected and every $s\in \Gamma(D\setminus
K,\,\mathcal{O})$ there exists a unique $\tilde{s}\in
\Gamma(D,\,\mathcal{O})$ such that $\tilde{s}=s$ on $D\setminus K$.
\end{proposition}

\medskip
\noindent ``Nice'' in the above statement means that the form $f$
can be smoothly extended on $U$ for some local embedding $\phi:
W\subset \Omega \to U^{\text{open}}\subset \mathbb{C}^N$.

\medskip
\noindent Hence to prove Theorem 1.3  it remains to show that we
can always find an open neighborhood $\Omega$ of $K$ that
satisfies property (P) of the above proposition. Since $K$ is
compact on an $(n-1)$-complete space we can always find an open,
relatively compact $(n-1)$-complete subdomain $\Omega$ of $X$ that
contains $K$. Property (P) for such domains has been established
by Theorem 4.2 (for $A:=\text{Sing}\,X$) and the proof of Theorem
1.3 is completed.

\section{A generalization of a lemma in \cite{FOV}}

\medskip
\noindent Proposition 2.2 in section 2 was a special case of the
following key lemma in \cite{FOV}:

\begin{lemma}(Lemma 2.1 in \cite{FOV}) For each $q>0$ and for each coherent,
torsion free $\mathcal{O}_{\tilde{X}}$-module $\mathcal{S}$ there
exists a\, $T\in \mathbb{N}$\, such that $i_{\tilde{\Omega},*}:
H^q(\tilde{\Omega}, J^T\mathcal{S})\to H^q(\tilde{\Omega},
\mathcal{S})$ is the zero map, where $i:
J^T\mathcal{S}\hookrightarrow \mathcal{S}$ is the inclusion map.
\end{lemma}

\medskip
\noindent In this section we shall see that it is possible to drop
the assumption that $S$ is torsion free and $J$\; is a principal
ideal sheaf in the lemma above. Hence lemma 6.1 should be valid for
any proper modification $\pi:\widetilde{X}\to X$ where is  $X$ a
pure $n$-dimensional reduced Stein space, $\Omega$ open relatively
compact Stein subdomain of
$X$,\;$\widetilde{\Omega}=:\pi^{-1}(\Omega)$,\; $S$ is any coherent
$\mathcal{O}_{\widetilde{X}}$-module and where $J$ is the ideal
sheaf of the exceptional set $\tilde{A}$ of the proper modification
$\pi$.

\medskip
\noindent{\bf{Remark:}} Recall that a proper modification between
two complex spaces $\widetilde{X},\,X$ consists of a proper
surjective holomorphic map $\pi: \widetilde{X}\to X$, closed
analytic subsets $\widetilde{A},\,A$ of $\widetilde{X},\,X$
respectively, such that a)\, $A=\pi(\widetilde{A})$,\;b)\,
$\widetilde{A},\,A$ are analytically rare, \;c)\,
$\pi_{\upharpoonright{\widetilde{X}\setminus
\widetilde{A}}}:\widetilde{X}\setminus \widetilde{A}\to X\setminus
A$ is a biholomorphism and d)\,$\widetilde{A},\,A$ are minimal with
respect to properties a)-c). For reduced spaces $X$ saying that a
closed analytic set $A$ is analytically rare is equivalent to saying
that no connected component of $X$ is contained in $A$. We would be
primarily interested in proper modifications $\pi: \widetilde{X}\to
X$, where $X$ is a reduced, pure $n$-dimensional Stein space. Then
$\widetilde{X}$ would also be reduced and pure $n$-dimensional (see
Chapter VII, page 287-288 in \cite{GPR}). If $\Omega$ is an open
relatively compact Stein subdomain of $X$, then we know that
$\Omega$ does not contain any compact $n$-dimensional irreducible
component. Hence $\widetilde{\Omega}$ can not contain any compact
$n$-dimensional irreducible component as well.

\medskip
\noindent The new ingredient in the proof is based on an application
of the Artin-Rees lemma:

\begin{lemma}\; (Artin-Rees Lemma, page 441 in \cite{J}) Let $M$ be a finitely generated
module over a noetherian ring $R$, and let $U$ be a  submodule of
$M$ and $I$ an ideal of $R$. Then there exists a positive integer
$k$ such that for all integers $n\ge k$ we have

$$
I^{n}\,M\cap U=I^{n-k}\,(I^k\,M\cap U).
$$
\end{lemma}

\medskip
\noindent For completion we shall recall the proof of lemma 6.1 and
mention all the necessary modifications.

\medskip
\noindent {\it{Proof of the more general version of Lemma 6.1.}} We
shall prove lemma 6.1 (under the more general assumptions on
$\pi,\,\mathcal{S},\,J$ mentioned above) using downward induction on
$q>0$. Observe that $\tilde{\Omega}$ is a pure $n$-dimensional
complex space with no compact $n$-dimensional branches. It follows
from the Main Theorem in Siu \cite{Siu} that $H^{n}(\tilde{\Omega},
\mathcal{S})=0$ for every coherent $\mathcal{O}_{\tilde{X}}$-module
$\mathcal{S}$. Hence, the statement is true for $q=n$ and any $T\in
\mathbb{N}$.

\medskip
\noindent  When $q>0,\,\,{\mbox{Supp}} R^q\pi_{*} \mathcal{S}$ is
contained in $A$. The annihilator ideal $\mathcal{A'}$ of
$R^q\pi_{*} \mathcal{S}$ is coherent and by Cartan's Theorem A there
exist functions \; $f_{1},\cdots, f_{L}\in \mathcal{A'}(X)$ that
generate each stalk $\mathcal{A'}_{z}$ in a neighborhood of
$\overline{\Omega}$. Let $\mathcal{A}$ be the
$\mathcal{O}_{\tilde{X}}$-ideal generated by $\tilde f_j=f_j\circ
\pi,\;\; 1\le j\le L$.  A crucial observation which will be useful
later, is that $(\tilde{f_j})_{\tilde{\Omega},*}:
H^q(\tilde{\Omega}, \mathcal{S}_{|_{\tilde{\Omega}}})\to
H^q(\tilde{\Omega}, \mathcal{S}_{|_{\tilde{\Omega}}})$ are zero for
all $j,\,1\le j\le L, \, q>0$. To see this, consider the following
commutative diagram

\[\begin{CD}
H^q(\tilde{\Omega}, \mathcal{S}_{|_{\tilde{\Omega}}}) @>(\tilde{f_j})_{\tilde{\Omega}, *}>>   H^q(\tilde{\Omega}, \mathcal{S}_{|_{\tilde{\Omega}}})\\
@V{\cong}VV                                                                   @V{\cong}VV\\
R^q\pi_*\mathcal{S}(\Omega)@>(f_j)_{\Omega,\#}>> R^q\pi_*\mathcal{S}(\Omega) \\
\end{CD}\]

\medskip
\noindent The vertical maps are isomorphisms, due to Satz 5, Section
2, in \cite{Gr}. Recalling the way $\mathcal{O}_X$ acts on
$R^q\pi_*\mathcal{S}$ and using the fact that the $f_j$'s are in the
annihilator ideal of $R^q\pi_*\mathcal{S}$ we conclude that
$(f_j)_{\Omega,\#}=0$. Hence, due to the commutativity of the above
diagram $(\tilde{f_j})_{\tilde{\Omega},*}$ is zero.

\medskip
\noindent Let $Z(\mathcal{A})\,({\rm resp.}\; Z(\mathcal{A'}))$
denote the zero variety of $\mathcal{A}\,({\rm resp.}\;
\mathcal{A'})$. Since
$Z(\mathcal{A'})={\mbox{Supp}}R^q\pi_*\mathcal{S}$ is contained in
$A,$ we have that $Z(\mathcal{A})$ is contained in $\tilde{A}$ near
$\overline{\tilde{\Omega}}$. Thus by R\"uckert's Nullstellensatz for
ideal sheaves, (see Theorem,  page 82 in \cite{GRem}), we have
$J^{\mu} \subset\mathcal{A}$ on $\tilde{\Omega}$  for some $\mu\in
\mathbb{N}$.
Consider the surjection $\phi: \mathcal{S}^{\oplus L} \to
\mathcal{A}\, \mathcal{S}$ given by $(s_1,\cdots,s_L)\to
\sum_{1}^{L} \tilde{f_{j}} s_j$ and set $K=\text{ker}\phi$. By
definition the sequence

\begin{eqnarray}\label{eq:shortexact}
0\to K \stackrel{i}\to \mathcal{S}^{\oplus L}\stackrel{\phi}\to
\mathcal{A}\, \mathcal{S}\to 0
\end{eqnarray}

\noindent is exact. Now, unlike the proof of Lemma 2.1  in
\cite{FOV} the following sequence will no longer be exact

$$
0\to J^a\, K \stackrel{i}\to J^a \, \mathcal{S}^{\oplus L}
\stackrel{\phi}\to J^a  \mathcal{A}\, \mathcal{S}\to 0.
$$

\noindent Recall that in our situation $S$ is not torsion free and
$J$ is not locally generated by one element. On the other hand for
every $a\ge 0$ the morphism $\phi: J^a\,\mathcal{S}^{\oplus\,L}\to
J^a\,\mathcal{A}\,\mathcal{S}$ is still surjective. Let
$K_a:=\text{kern}(\phi)_{|_{J^a\,\mathcal{S}^{\oplus\,L}}}=K\cap
J^a\,\mathcal{S}^{\oplus\,L}$. Then the following sequence will be
exact

$$
0\to K_a\overset{i}\to J^a\,\mathcal{S}^{\oplus L}\overset{\phi}\to
J^a\mathcal{A}\,\mathcal{S}\to 0.
$$

\medskip
\noindent Taking all the above into consideration we obtain the
following commutative diagram:



$$
\xymatrix{
 & H^q(\tilde{\Omega}, J^{a+\mu}\mathcal{S}) \ar[d] \\
 & H^q(\tilde{\Omega}, J^a \mathcal{A}\,
\mathcal{S})\ar[d] \ar[r]^{\delta}  & H^{q+1}(\tilde{\Omega},
 K_a)\ar[d]_{i_1}\\
H^q(\tilde{\Omega}, \mathcal{S})^{\oplus L}\ar[dr]^{\chi}
\ar[r]^{\phi_{\tilde{\Omega}, *}}  & H^q(\tilde{\Omega},
\mathcal{A}\, \mathcal{S}) \ar[d]^{i_2}
\ar[r]^{\delta}  & H^{q+1}(\tilde{\Omega}, K)\\
& H^q(\tilde{\Omega}, \mathcal{S}) }
$$

\bigskip
\noindent where the third row  is exact (as part of the long exact
cohomology sequence that arises from (\ref{eq:shortexact})\,) and
the vertical maps are induced by sheaf inclusions. The map $\chi$ is
defined to be $\chi:=i_2\circ \phi_{\tilde{\Omega},\,*}$ and we can
show that $\chi(c_1,\cdots,c_L)=\sum_{j=1}^{L}
(\tilde{f_j})_{\tilde{\Omega},\,*} c_j,$ where\; $c_j\in
H^q(\tilde{\Omega}, \mathcal{S}),\; 1\le j\le L$.

\smallskip
\noindent In the original commutative diagram in \cite{FOV} the
cohomology group in the second line-third column was
$H^{q+1}(\tilde{\Omega},\,J^a\,K)$. The induction hypothesis applied
to $K$, allowed us to claim that there exists an $a$ large enough so
that $i_1=0$. Here, we {\it{need}} to show that $i_1:
H^{q+1}(\tilde{\Omega},\,K_a)\to H^{q+1}(\tilde{\Omega},\,K)$ is the
zero map assuming by the induction hypothesis that the map $H^{q+1}
(\tilde{\Omega},\,J^{\alpha}\,K)\overset{I}\to
H^{q+1}(\tilde{\Omega},\,K)$ is the zero map for some $\alpha$
sufficiently large. If we were able to show that the map $i_1$
factors through $H^{q+1}(\tilde{\Omega},\,J^a\,K)$ then we would be
done since the rest of the proof carries over verbatim.

\medskip
\noindent For each positive integer $k$ we set

$$A_k:=\text{supp}\,\left(\frac{(J^{\alpha+k}\mathcal{S}^{\oplus L}\cap
K)}{J^{\alpha}\,(J^k\,\mathcal{S}^{\oplus L}\cap K)}\right).
$$

\smallskip
\noindent Each $A_k$ is a subvariety of $\widetilde{X}$ and by
Artin-Rees lemma we have that $\cap_{k=1}^{\infty} A_k=\emptyset$.
By the descending chain property of analytic varieties,  there
exists an integer $m$ such that $\tilde{\Omega}\cap
\left(\cap_{k=1}^{m} A_k\right)=\emptyset$, since
$\tilde{\Omega}\Subset \tilde{X}$. Hence for every $x\in
\tilde{\Omega}$ there exists $k_x\le m$ such that

$$
K_{k_x+\alpha,\,x}=J_x^{k_x+\alpha}\,\mathcal{S}_x^{\oplus L}\cap
K_x=J^{\alpha}_x\,(J^{k_x}_x\,\mathcal{S}^{\oplus L}_x \cap
K_x)\subset J^{\alpha}_x\, K_x.
$$

\medskip
\noindent Let us choose $a:=m+\alpha$. Then we have
$K_{a,\,x}\subset K_{k_x+\alpha,\,x}\subset J^a_x\,K_x$. Thus over
$\tilde{\Omega}$, the inclusion $K_{a}\to K$ factors via
$J^{\alpha}\,K$ so the map $i_1: H^{q+1}(\tilde{\Omega},\,K_a)\to
H^{q+1}(\widetilde{\Omega},\,K)$ is the zero map.

\smallskip
\noindent Then, for an element $\sigma \in H^q(\tilde{\Omega},
\mathcal{A}\, \mathcal{S})$ that comes from $H^q(\tilde{\Omega},
J^{a+\mu} \mathcal{S})$, we have $\delta\sigma=0$, so
$\sigma=\phi_{\tilde{\Omega},*}(\sigma_1,\cdots, \sigma_L)$  for
\;$\sigma_j\in H^q(\tilde{\Omega}, \mathcal{S}),\,1\le j\le L$. By
the crucial observation above and the way $\chi$ is defined, we
conclude that $\chi$ is the zero map. Hence $i_2(\sigma)=i_2\circ
\phi_{\tilde{\Omega},*}(\sigma_1,\cdots, \sigma_L)=\sum_{j=1}^{L}
(\tilde{f_j})_{\tilde{\Omega},*} \sigma_j=0$. Thus, for $i:J^{a+\mu}
\mathcal{S}\hookrightarrow \mathcal{S}$ the inclusion map, we have
that $i_{\tilde{\Omega},*}: H^q(\tilde{\Omega},
J^{a+\mu}\mathcal{S})\to H^q(\tilde{\Omega}, \mathcal{S})$ is the
zero map.
\medskip
\noindent

\medskip
\noindent

\section{Appendix}

\medskip
\noindent  The purpose of this appendix is to prove the equivalence
of the corresponding distance functions $d_A$ over $K$ that was
alluded to in section 2.1., if we choose different compacts
$K'\Supset K$ and different embeddings.

\begin{lemma} Let $\phi: U\to \mathbb{C}^T$ be a local embedding and
let $K,\,L,\,M$ be compact subsets of $U$ such that $K\Subset
L\Subset U$  and $K\Subset M \Subset U$. Let
$d_1(x):=\text{dist}(\phi(x),\,\phi(L\cap A))$ and
$d_2(x):=\text{dist}(\phi(x),\,\phi(M\cap A))$. Then

$$
d_1(x)\sim_K d_2(x).
$$
\end{lemma}

\medskip
\noindent{\it{Proof:}} If $M\subset L$ then clearly $d_1(x)\le
d_2(x)$. Suppose that $M\not\subset L$. Then $\text{dist}(\phi(K),\,
\phi(M\setminus \overset{0} L\cap A))=a>0$ (here $\overset{0} L$
denotes the interior of the set $L$). Let us assume that there
existed $x\in K$ such that $d_2(x)<d_1(x)$. Then we must have that
$d_2(x)\ge a$, hence $1\le a^{-1}\,d_2(x)$. Let $B:=\text{max}_{x\in
K} d_1(x)$. Then we have

$$
d_1(x)\le B=B\times 1\le B\,a^{-1} d_2(x).$$

\medskip
\noindent Choosing as $C_1:=\text{max}(\{1,\,B\,a^{-1}\})$ we obtain
for all $x\in K$ that $d_1(x)\le C_1\,d_2(x)$. Similarly we can
prove that there exists  a constant $C_2\ge 1$  such that $d_2(x)\le
C_2\,d_1(x)$ for all $x\in K$. Choosing as
$C:=\text{max}(C_1,\,C_2)$ we obtain the equivalence of $d_1$ and
$d_2$ over $K$.

\medskip
\noindent
\begin{lemma} Let $\phi_i: U_i\to \mathbb{C}^{T_i}$ be local
embeddings for $i=1,\,2$ and let $K,\,K'$ be a compact subsets of
$U_1\cap U_2$ such that $K\Subset K'\Subset U_1\cap U_2$. Then we
have:

$$
\text{dist}(\phi_1(x),\,\phi_1(K'\cap A))\sim_K
\text{dist}(\phi_2(x),\,\phi_2(K'\cap A)).
$$
\end{lemma}

\medskip
\noindent{\it{Proof:}} Let $x_0\in K$ and $(\phi_0,\,V_0, W_0)$ be a
minimal embedding of $X$ at $x_0$; $\phi_0: V_0\to W_0\subset
\mathbb{C}^{N_0}$. In this case we know that $d\,\phi_{0,\,x_0}:
T_{x_0}X\to T_{\phi_0(x_0)} W_0$ is an isomorphism. By shrinking
$V_0$ we can assume that $V_0\subset U_1\cap U_2$. Shrinking further
$V_0,\,W_0$ and after a holomorphic change of coordinates we can
assume that $\phi_0(x_0)=0=\phi_i(x_0)$ for $i=1,\,2$.

\medskip
\noindent {\bf{Claim:}} The maps $\phi_i\circ
\phi^{-1}_0:\phi_0(V_0)\to \phi_i(V_0)$ for $i=1,\,2$ extend from
$\phi_0(V_0)$ to an embedding $\tilde{\phi_i}: \mathbb{C}^{N_0}\to
\mathbb{C}^{T_i}$.

\medskip
\noindent {\it{Proof of Claim:}} Indeed, as $\phi_i$ is an embedding
of $X$ at $x_0$ we know that $d\,\phi_{i,\,x_0}: T_{x_0}X\to
T_{0}\mathbb{C}^{T_i}$ is injective. Let $\phi_i\circ \phi^{-1}_0:
\phi_0(V_0)\to \phi_i(V_0)$ be given by the coordinate map
$(f_1,\,\cdots,f_{T_i})$. Then $\phi_i\circ \phi^{-1}_0$ extends to
a map $\tilde{\phi_i}: \mathbb{C}^{N_0}\to \mathbb{C}^{T_i}$ by
extending each $f_j$ for $j=1,\cdots, T_i$. Then we have the
following commutative diagram

$$
\xymatrix{ & T_{x_0}X
 \ar[d]_{d\,\phi_{0,\,x_0}}^{\cong} \ar[r]^{d\phi_{i,\,x_0}} &T_0\mathbb{C}^{T_i}\\
&T_{0}\mathbb{C}^{N_0} \ar[ur]_{d\tilde{\phi}_{i,\,0}}
 }
$$

\medskip
\noindent Taking into account the construction of the
$\tilde{\phi_i}$ the fact that $d\phi_0$ is an isomorphism at $x_0$
and each $d\phi_i$ is an embedding of $X$ at $x_0$ we obtain that
$d\tilde{\phi_i}_{x_0}$ is injective; hence $\tilde{\phi_i}$ is an
embedding at $x_0$ for $i=1,\,2$.

\medskip
\noindent When $L_0\Subset W_0$ is a neighborhood of $\phi_0(x_0)=0$
it follows that

$$
|\tilde{\phi_i}(x)-\tilde{\phi_i}(y)|\sim_K |x-y|
$$

\noindent on $L_0\times L_0$ and for $i=1,\,2$. Hence, when $K_0$ is
a compact neighborhood of $x_0$ contained in
${\overset{0}{K'}}_0:={\phi_0}^{-1}({\overset{0}L}_0)\;\;
({\overset{0}L}_0\;\text{or}\; {\overset{0}{K'}}_0$ denotes the
interior of the corresponding sets) then we have

$$
\text{dist}(\phi_i(x),\,\phi_i(K'_0\cap
A))\sim_{K_0}\;\text{dist}(\phi_0(x),\,\phi_0(K'_0\cap A)).
$$

\medskip
\noindent Using the transitivity of $\sim_{K_0}$ and Lemma 7.1 we
obtain that

$$
\text{dist}(\phi_1(x),\,\phi_1(K'\cap
A))\sim_{K_0}\,\text{dist}(\phi_2(x),\,\phi_2(K'\cap A)).
$$

\medskip
\noindent We can now cover $K$ by finitely many $K_0$ and obtain the
desired result.

\end{document}